\documentclass[11pt]{article}
\usepackage[margin=1in]{geometry}
\usepackage[T2A,T1]{fontenc}
\usepackage[russian,english]{babel}
\usepackage{lmodern}
\usepackage{amsmath,amssymb,amsthm,microtype}
\usepackage{tikz}
\usepackage[hidelinks]{hyperref}
\newtheorem{theorem}{Theorem}
\newtheorem{lemma}{Lemma}
\title{Three Circles in a Unit Disk: Three Solutions to an Extremal Problem}
\author{Anatoly Eydelzon\\
\small Department of Mathematical Sciences\\
\small The University of Texas at Dallas\\
\small Email: \href{mailto:anatoly@utdallas.edu}{anatoly@utdallas.edu}\\
\small ORCID: \href{https://orcid.org/0000-0002-8822-9277}{0000-0002-8822-9277}}
\date{November 2026}
\begin{document}
\maketitle
\begin{abstract}
We place two disks inside the unit disk so that they touch each other and the unit circle. We then insert a disk of the largest possible radius, requiring it to touch both original disks. We show that the smallest possible sum of the three areas is $11\pi/18$, and find all equality cases. We give three solutions, using Descartes' circle theorem, inversion, and contact-point coordinates. The main difficulty is checking that the third disk is largest among all disks that fit, not only among disks touching the first two.
\end{abstract}

\section*{Publication Status}
\begin{otherlanguage}{russian}
Русскоязычная версия настоящей работы представлена в журнал
«Математическое просвещение» и находится на рассмотрении редакции.
\end{otherlanguage}
\begin{theorem}\label{thm:main}
Let two disks of positive radii be contained in the unit disk, tangent to each other and internally tangent to the unit circle. Insert a third disk of \emph{largest possible radius} without overlapping either of the first two, and suppose this largest disk also touches both original disks. Then the sum of the areas of the three disks is at least
\[
{\frac{11\pi}{18}}.
\]
Equality holds precisely when the three radii are $\{1/3,1/2,1/2\}$. There are two ways to obtain equality: the original radii are $(1/2,1/2)$ and the inserted radius is $1/3$, or the original radii are $(1/3,1/2)$ and the inserted radius is $1/2$, up to interchanging the original disks and rigid motions.
\end{theorem}

At first glance, one might expect the minimum sum of areas to occur when all three disks have the same radius. Surprisingly, this is not the case. Despite the symmetry of the problem, the minimizing configuration is asymmetric: the radii are $1/2$, $1/2$, and $1/3$, and the minimum total area is $11\pi/18$. Moreover, the same minimum is achieved with two different orders of insertion. Thus the simplicity of the final answer contrasts with the less obvious geometry behind it.

\begin{figure}[htbp]
\centering
\begin{tikzpicture}[scale=2.05, line cap=round, line join=round]
  \begin{scope}[shift={(-1.42,0)}]
    \draw[thick] (0,0) circle (1);
    \fill[blue!10] (0.628362,-0.351513) circle (0.28);
    \draw[blue!65!black,thick] (0.628362,-0.351513) circle (0.28);
    \fill[orange!14] (0.541089,0.302692) circle (0.38);
    \draw[orange!75!black,thick] (0.541089,0.302692) circle (0.38);
    \fill[green!12] (-0.263241,-0.281877) circle (0.614318);
    \draw[green!45!black,very thick] (-0.263241,-0.281877) circle (0.614318);
    \draw[gray!75,densely dashed] (0.884538,-0.052271) circle (0.113919);
    \fill (0,0) circle (0.016) node[above left,scale=.63] {$O$};
    \fill (0.628362,-0.351513) circle (0.015) node[below right,scale=.67] {$A$};
    \fill (0.541089,0.302692) circle (0.015) node[above right,scale=.67] {$B$};
    \fill (-0.263241,-0.281877) circle (0.015) node[below left,scale=.67] {$C$};
    \node[scale=.78] at (0,-1.28) {(a) A typical configuration};
  \end{scope}
  \begin{scope}[shift={(1.42,0)}]
    \draw[thick] (0,0) circle (1);
    \fill[blue!10] (0,-0.5) circle (0.5);
    \draw[blue!65!black,thick] (0,-0.5) circle (0.5);
    \fill[orange!14] (0,0.5) circle (0.5);
    \draw[orange!75!black,thick] (0,0.5) circle (0.5);
    \fill[green!12] (0.666667,0) circle (0.333333);
    \draw[green!45!black,very thick] (0.666667,0) circle (0.333333);
    \fill (0,0) circle (0.016) node[above left,scale=.63] {$O$};
    \fill (0,-0.5) circle (0.015) node[left,scale=.67] {$A$};
    \fill (0,0.5) circle (0.015) node[left,scale=.67] {$B$};
    \fill (0.666667,0) circle (0.015) node[above,scale=.67] {$C$};
    \node[scale=.78] at (0,-1.28) {(b) Equality: $a=b=\tfrac12$, $c=\tfrac13$};
  \end{scope}
\end{tikzpicture}
\caption{The original disks (blue and orange) and the largest disk that fits (green). The dashed circle in (a) is the other circle tangent to the two original disks and the unit circle. The figure is drawn to scale.}
\label{fig:configurations}
\end{figure}
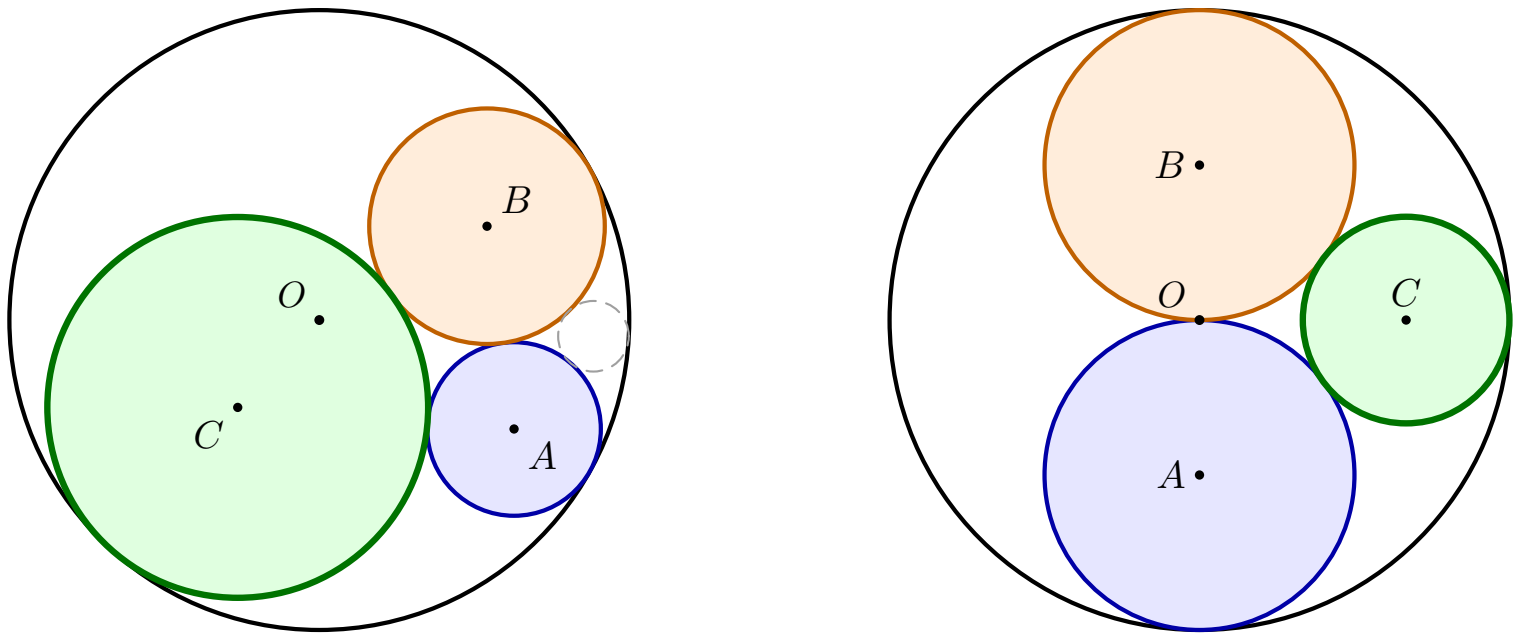

Figure~\ref{fig:configurations} shows a typical arrangement and an equality case. In the left panel, the dashed circle touches both original disks and the unit circle, but is smaller than the disk we insert.

\paragraph{Background.}
The relation for four mutually tangent circles goes back to Descartes \cite{descartes} and was rediscovered by Soddy \cite{soddy}; see \cite{lagarias,graham} for modern treatments, in which the enclosing circle carries negative curvature. The second solution uses classical circle inversion \cite{coxeter,johnson}. The third uses contact points to describe a family of tangent circles. To the best of our knowledge, the precise sequential extremal problem considered here has not previously appeared in the literature.

\section{A global maximality criterion}
In this section and the first solution, let $a$ and $b$ denote the radii of the two original disks, labeled so that $0<a\le b$. Since the disks are internally tangent to the unit circle and externally tangent to each other, $a+b\le 1$. Let $c$ denote the radius of the third, inserted disk.
\begin{lemma}[Global maximality and the admissible region]\label{lem:global}
Put
\[
L(b)=\frac{b(1-b)}{1-b+b^2}.
\]
The largest possible third radius is
\begin{equation}\label{eq:cmax}
c_{\max}=\begin{cases}
1-b,&a\le L(b),\\[2mm]
\displaystyle\frac{ab}{a+b-ab-2\sqrt{ab(1-a-b)}},&a\ge L(b).
\end{cases}
\end{equation}
The two expressions agree at $a=L(b)$. A largest third circle is tangent to both original circles and the unit circle if and only if $a\ge L(b)$. Thus the parameter region in the theorem is
\begin{equation}\label{eq:region}
0<b<1,\qquad L(b)\le a\le\min\{b,1-b\}.
\end{equation}
\end{lemma}
\begin{proof}
Let $O$ be the unit disk's center and let $A,B$ be the original centers, with $OA=1-a$, $OB=1-b$, and $AB=a+b$. The smaller angle $\theta=\angle AOB\in(0,\pi]$ satisfies
\begin{equation}\label{eq:sin}
\sin^2\frac\theta2=\frac{ab}{(1-a)(1-b)}.
\end{equation}
Any feasible disk of radius $c$ satisfies $c\le1-b$, since its center $X$ obeys
\[
b+c\le |X-B|\le |X|+|B|\le(1-c)+(1-b).
\]
For fixed $0<c\le1-b$, it suffices to test centers on $|X|=1-c$. Indeed, the forbidden open center-disks $D(A,a+c)$ and $D(B,b+c)$ overlap, since $AB=a+b<(a+c)+(b+c)$. We can find a path to infinity that avoids both disks. Choose a point $P$ in their intersection (which is nonempty by the strict inequality above). For any feasible center $X$ outside both forbidden open disks, consider the ray
\[
Y(s)=P+s(X-P),\qquad s\ge1.
\]
Each forbidden disk is convex and contains $P$. Its intersection with the ray starting at $P$ is therefore an initial interval: if $Y(s)$ belonged to that disk for some $s>1$, convexity would imply that $X=Y(1)$ also belonged to it, a contradiction. Hence $Y(s)$ stays outside both forbidden open disks for every $s\ge1$, and $|Y(s)|\to\infty$. If $|X|<1-c$, continuity gives a first parameter $s_*\ge1$ such that $|Y(s_*)|=1-c$. The segment from $X=Y(1)$ to $Y(s_*)$ stays inside the closed disk $|Y|\le1-c$: both endpoints lie in this disk, which is convex. As shown above, the entire segment also avoids both forbidden open disks. Thus the trial center remains feasible until it reaches the outer boundary. We have therefore found a disk of the same radius that also touches the unit circle. So, for each fixed radius, we need only check disks that touch the outer circle.

For an outer-tangent trial disk of radius $c$, the disk of radius $r\in\{a,b\}$ forbids an open angular interval about its own outer contact direction, of half-width
\begin{equation}\label{eq:alpha}
\alpha_r(c)=2\arcsin\sqrt{\frac{rc}{(1-r)(1-c)}}.
\end{equation}
This follows by applying the cosine rule to the two centers and the origin. The two contact directions are separated by arcs of lengths $\theta$ and $2\pi-\theta$; see Figure~\ref{fig:configurations}(a). Thus a trial disk fits if and only if the forbidden intervals leave a point on at least one arc, equivalently
\begin{equation}\label{eq:fit}
\alpha_a(c)+\alpha_b(c)\le2\pi-\theta.
\end{equation}
Indeed $\theta\le\pi$, so the longer arc is decisive. Both half-widths increase strictly with $c$ (up to their endpoint values).

At $c=1-b$, the only possible center is antipodal to $B$. Choose $B=(1-b,0)$ and $A=(1-a)(\cos\theta,\sin\theta)$; then the candidate center is $X=(-b,0)$. Using \eqref{eq:sin},
\begin{equation}\label{eq:anti}
|X-A|^2-(a+1-b)^2
=\frac{4\bigl[b(1-b)-a(1-b+b^2)\bigr]}{1-b}.
\end{equation}
Therefore this antipodal disk fits precisely when $a\le L(b)$. If $a>L(b)$, strict monotonicity in \eqref{eq:fit} gives a unique limiting $c<1-b$ for which the two forbidden intervals meet. The limiting circle touches both original circles and the unit circle, with its outer contact point on the longer arc. To identify the Descartes root, observe that a circle tangent to all three has its outer contact point on either the shorter or the longer arc. The shorter-arc contact gives the equation $\alpha_a(c)+\alpha_b(c)=\theta$, whereas the longer-arc contact gives $\alpha_a(c)+\alpha_b(c)=2\pi-\theta$. Since the left-hand side is strictly increasing in $c$ and $\theta\le 2\pi-\theta$, the longer-arc circle has the larger radius (the two radii coincide if $\theta=\pi$). Solving the Descartes relation (see, e.g., \cite{soddy,lagarias}) gives the two roots $\frac{ab}{a+b-ab\mp2\sqrt{ab(1-a-b)}}$; the larger is the second expression in \eqref{eq:cmax}. At $a=L(b)$ it equals $1-b$. If $a<L(b)$ the maximal disk touches only the $b$-disk and the unit circle, so does not meet the theorem's three-contact requirement. This proves \eqref{eq:cmax}--\eqref{eq:region}.
\end{proof}
\section{First solution: the Descartes tangency surface}
We first use Descartes' circle theorem to turn the problem into an algebraic minimization.
For the theorem, on the admissible region \eqref{eq:region}, write
\begin{equation}\label{eq:F}
F(a,b)=a^2+b^2+c(a,b)^2,\qquad
c(a,b)=\frac{ab}{a+b-ab-2\sqrt{ab(1-a-b)}}.
\end{equation}
The denominator is positive for $a,b>0$ and $a+b\le1$. We prove $F\ge11/18$.

\paragraph{Stationary points on the Descartes surface.}
The following identity is the Descartes relation for four mutually tangent circles, with the enclosing unit circle assigned oriented curvature $-1$ (see \cite{lagarias,graham}). On the smooth Descartes surface define
\begin{equation}\label{eq:H}
H(a,b,c)=(ab+bc+ca+abc)^2-4abc(a+b+c)=0.
\end{equation}
For $a,b,c>0$ this surface is smooth at every configuration with $a+b<1$ and distinct Descartes roots: indeed the two roots in $c$ are distinct, so $H_c\ne0$. At a constrained stationary point of $a^2+b^2+c^2$ we have $2a=\lambda H_a$, $2b=\lambda H_b$, $2c=\lambda H_c$. If $a,b,c$ were all distinct, subtraction and division would give
\[
\frac{H_a-H_b}{a-b}=\frac{H_b-H_c}{b-c}.
\]
But direct expansion gives
\begin{equation}\label{eq:diff}
\frac{H_a-H_b}{a-b}-\frac{H_b-H_c}{b-c}
=2(a-c)\bigl(ab+bc+ca+abc-2(a+b+c)\bigr).
\end{equation}
The last parenthesis is strictly negative when $0<a,b,c<1$: $ab+bc+ca<a+b+c$ and $abc<a+b+c$. Thus every regular stationary point has at least two equal radii.

To find stationary points, we work on the whole positive Descartes surface, not just the part allowed by our problem. This may give extra candidates, but it does not miss any stationary point in the allowed region. By symmetry of \eqref{eq:H} and the objective, it suffices to study $a=b=r$. Set $r=(1-t^2)/2$, $0<t<1$. The two Descartes radii are
\[
c_+(t)=\frac{1+t}{3-t},\qquad c_-(t)=\frac{1-t}{3+t}.
\]
For $F_\pm(t)=\frac12(1-t^2)^2+c_\pm(t)^2$, direct differentiation yields
\begin{align*}
F_+'(t)&=\frac{2(t+1)(t^2-4t+1)(t^3-6t^2+11t-4)}{(t-3)^3},\\
F_-'(t)&=\frac{2(t-1)(t^2+4t+1)(t^3+6t^2+11t+4)}{(t+3)^3}.
\end{align*}
The minus branch has no stationary point for $0<t<1$. The plus branch has precisely two: $t=2-\sqrt3$ and the unique zero $t_*$ of $t^3-6t^2+11t-4$ in $(0,1)$ (the cubic is strictly increasing there, since its derivative is $3t^2-12t+11>0$). The corresponding objective values are $63-36\sqrt3>11/18$ and $F_+(t_*)>11/18$, respectively. For the second value, direct substitution into $F_+$ gives
\[
F_+(t_*)-\frac{11}{18}=\frac{t_*K(t_*)}{18(3-t_*)^2}>0,
\]
where the polynomial $K$ is displayed in \eqref{eq:K} and proved strictly positive on $[0,1)$ immediately below. Thus this comparison does not rely on a numerical approximation.

\paragraph{Boundary $a=L(b)$.}
Here $c=1-b$ and
\begin{equation}\label{eq:Q}
F-\frac{11}{18}
=\frac{(2b-1)^2Q(b)}{18(1-b+b^2)^2},\qquad
Q(b)=9\left(b-\frac12\right)^4+\frac{35}{2}\left(b-\frac12\right)^2+\frac{33}{16}>0.
\end{equation}
Thus equality on this boundary occurs only at $b=1/2$, giving $(a,b,c)=(1/3,1/2,1/2)$; here the third disk is as large as the larger initial disk.

\paragraph{Boundary $a=b$ (necessarily $b\le1/2$).}
With $t=\sqrt{1-2b}\in[0,1)$ and $c=c_+(t)$, one has
\begin{equation}\label{eq:K}
F-\frac{11}{18}=\frac{t\,K(t)}{18(3-t)^2},\quad
K(t)=9t^5-54t^4+63t^3+108t^2-146t+48.
\end{equation}
To show $K(t)>0$ for $0\le t<1$, put $t=z/(1+z)$, $z\ge0$. Then
\[
(1+z)^5K\left(\frac{z}{1+z}\right)
=28z^5+52z^4-9z^3+4z^2+94z+48>0,
\]
since $52z^4+4z^2\ge 2\sqrt{208}\,z^3>9z^3$ for $z>0$, while the value at $z=0$ is $48$. Equality in \eqref{eq:K} occurs only at $t=0$, i.e., $(a,b,c)=(1/2,1/2,1/3)$; here the third disk is smaller than both initial disks. Thus the two equality cases have identical unordered radii but genuinely different insertion orders.

\paragraph{Boundary $a+b=1$ (necessarily $b\ge1/2$).}
Here the two Descartes roots coincide and $c=ab/(1-ab)$. With $u=ab\in(0,1/4]$,
\[
F=1-2u+\frac{u^2}{(1-u)^2},\qquad
\frac{dF}{du}=-2+\frac{2u}{(1-u)^3}<0
\]
for $0<u\le1/4$. Thus the minimum is at $u=1/4$, yielding $(a,b,c)=(1/2,1/2,1/3)$.

\paragraph{Endpoints and global conclusion.}
The region \eqref{eq:region} is bounded but has two excluded limiting endpoints: $b\to0$ and $b\to1$. As $b\to1$, $F\ge b^2\to1$. As $b\to0$, the inequalities $L(b)\le a\le b$ imply $a/b\to1$. Let $k=a/b$. Rationalizing \eqref{eq:F} gives
\[
c=\frac{kb\bigl(k+1-kb+2\sqrt{k(1-(k+1)b)}\bigr)}{(k-1)^2+2kb(k+1)+k^2b^2}.
\]
Since $L(b)/b\le k\le1$, we have $1-k=O(b^2)$. The numerator is asymptotic to $4b$ and the denominator to $4b$, so $c\to1$ and $F\to1$. These are the only missing endpoints of the bounded parameter region. At both endpoints, $F$ tends to $1>11/18$. Thus the fact that the region is not closed causes no difficulty: a minimum must occur at an interior stationary point or on one of the three boundary pieces considered above. Every stationary value exceeds $11/18$, and the boundary minimum is exactly $11/18$. Multiplying by $\pi$ proves the theorem:
\[
{\min\bigl(\pi a^2+\pi b^2+\pi c_{\max}^2\bigr)=\frac{11\pi}{18}}.
\]

\section{Second solution: inversion into a strip}
Our second proof uses inversion to replace two tangent circles by parallel lines (see \cite{coxeter,johnson}). This gives new coordinates in which we can minimize at the endpoints.
\subsection{Inversion at the north pole}
Rotate the picture so that the $a$-disk touches the unit circle at $N=(0,1)$. Now invert about the circle of radius $2$ centered at $N$:
\[
 (x,y)\longmapsto (X,Y)=(0,1)+\frac{4(x,y-1)}{x^2+(y-1)^2}.
\]
The unit circle becomes the line $Y=-1$, and the $a$-circle becomes $Y=1-2/a$. For the three-contact configuration under consideration, the images of the second original disk and the proposed third disk are congruent disks, each tangent to both parallel lines and to the other. The two inverted disks are congruent because they fit between the same parallel lines. But this alone does not show that the third disk is largest in the original configuration. We check that in Lemma~\ref{inv:lem:domain}. Their common inverted radius is
\[
 R=\frac{1-a}{a}>0.
\]
Their centers can be written $(u-R,-1-R)$ and $(u+R,-1-R)$, where the first corresponds to $b$ and the second to $c$. Inverting back gives
\begin{equation}\label{inv:eq:radii}
 a=\frac1{R+1},\qquad
 b=\frac{4R}{(u-R)^2+4(R+1)},\qquad
 c=\frac{4R}{(u+R)^2+4(R+1)}.
\end{equation}
Here we used the fact that a circle of radius $R$ centered at $(t,-1-R)$ inverts to a circle of radius $4R/[t^2+4(R+1)]$.

\begin{lemma}[Admissible inversion parameters]\label{inv:lem:domain}
The third disk in \eqref{inv:eq:radii} is a largest insertable disk, touching both originals, precisely when
\begin{equation}\label{inv:eq:udomain}
 -R\le u\le R,\qquad u^2\ge R^2-4.
\end{equation}
Consequently, with $v=u^2$, the admissible interval for each $R>0$ is
\begin{equation}\label{inv:eq:vdomain}
 \max\{0,R^2-4\}\le v\le R^2.
\end{equation}
Conversely, every $v$ in \eqref{inv:eq:vdomain} is realized by at least one admissible choice of sign of $u$.
\end{lemma}
\begin{proof}
In the inverted strip, the two circles tangent to both boundary lines
and to the transformed $b$-circle have horizontal center coordinates
$u+R$ and $u-3R$. Inverting back, their radii are respectively
\[
 c_+=\frac{4R}{(u+R)^2+4(R+1)},\qquad
 c_-=\frac{4R}{(u-3R)^2+4(R+1)}.
\]
Of these two candidates, the proposed third circle is the larger
exactly when $c_+\ge c_-$, or equivalently $u\le R$.

But we still have to compare it with \emph{every} possible
third disk, including disks that do not touch both originals. The geometric criterion of the preceding section is symmetric
in the original radii: regardless of their order, the larger
three-contact disk is globally largest precisely when
\begin{equation}\label{inv:eq:both-thresholds}
 b\ge L(a),\qquad a\ge L(b).
\end{equation}
To see this, suppose first that $a\le b$. The first inequality holds because
$L(a)<a\le b$, while the second is exactly the criterion already
proved. The case $b\le a$ is the same with the roles reversed.

Set $D_b=(u-R)^2+4(R+1)>0$. Substitution of \eqref{inv:eq:radii}
into \eqref{inv:eq:both-thresholds} gives the exact identities
\[
 b-L(a)=\frac{R(3R-u)(R+u)}{(R^2+R+1)D_b},
\]
and
\[
 a-L(b)=
 \frac{(u^2+4-R^2)((u-2R)^2+4-R^2)}
 {(R+1)(D_b^2-4RD_b+16R^2)}.
\]
The denominators are positive. In particular,
$D_b^2-4RD_b+16R^2=(D_b-2R)^2+12R^2>0$.
For $u\le R$, the factor $3R-u$ is positive, so
$b\ge L(a)$ is equivalent to $u\ge-R$. Once
$-R\le u\le R$, we have
\[
 (u-2R)^2+4-R^2\ge4>0,
\]
and therefore $a\ge L(b)$ is equivalent to
$u^2+4-R^2\ge0$. Combining the conditions gives exactly
\[
 -R\le u\le R,\qquad u^2\ge R^2-4.
\]
The argument also works in reverse, so the conditions are both
necessary and sufficient. Finally, $v=u^2$ gives
\eqref{inv:eq:vdomain}; every such $v$ is realized, for example, by
$u=-\sqrt v$.
\end{proof}

\subsection{An endpoint principle}
Put $v=u^2$ and $D=v+(R+2)^2$. The total area divided by $\pi$ is
\begin{equation}\label{inv:eq:F}
 F(R,v)=\frac1{(R+1)^2}+\frac{32R^2(D^2+4R^2v)}{(D^2-4R^2v)^2}.
\end{equation}
The denominator is positive, since $D^2-4R^2v$ is the product of the two positive denominators in \eqref{inv:eq:radii}.

\begin{lemma}[No interior minimum]\label{inv:lem:endpoint}
For fixed $R>0$, the minimum of $F(R,v)$ on the interval \eqref{inv:eq:vdomain} occurs at an endpoint.
\end{lemma}
\begin{proof}
Differentiating with respect to $v$ gives
\[
 \frac{\partial F}{\partial v}=-\frac{64R^2P_R(v)}{(D^2-4R^2v)^3},
\]
where
\begin{align*}
 P_R(v)={}&v^3+(9R^2+12R+12)v^2\\
 &+(-5R^4+24R^3+72R^2+96R+48)v\\
 &-5R^6-36R^5-84R^4-32R^3+144R^2+192R+64.
\end{align*}
Now $P_R''(v)=6v+18R^2+24R+24>0$. If $0<R\le2$, then
\[
 P_R'(0)=R^3(24-5R)+72R^2+96R+48>0.
\]
If $R>2$, then
\[
 P_R'(R^2-4)=16R^3(R+3)>0.
\]
Thus $P_R$ is strictly increasing on the admissible interval. Since $F_v$ has the opposite sign, it can change sign only from positive to negative. An interior critical point can therefore only be a maximum. The minimum must occur at one of the endpoints.
\end{proof}

\subsection{Endpoint evaluation and equality}
We now check the endpoints. There are three cases, depending on $R$.

\smallskip\noindent\emph{Upper endpoint $v=R^2$, $R>0$.} Factoring the difference gives
\begin{equation}\label{inv:eq:upper}
 F(R,R^2)-\frac{11}{18}=
 \frac{(R-1)^2(7R^4+6R^3+7R^2+6R+7)}{18(R+1)^2(R^2+R+1)^2}\ge0.
\end{equation}
Equality holds exactly when $R=1$.

\smallskip\noindent\emph{Lower endpoint $v=0$, $0<R\le2$.} Here
\begin{equation}\label{inv:eq:lower0}
 F(R,0)-\frac{11}{18}=\frac{(2-R)Q(R)}{18(R+1)^2(R+2)^4},
\end{equation}
where
\[
 Q(R)=11R^5+132R^4+121R^3-86R^2-36R+56.
\]
We need $Q(R)>0$ for $0<R\le2$. Set $R=2z/(1+z)$, $z>0$. Then
\begin{align*}
 (1+z)^5Q\!\left(\frac{2z}{1+z}\right)
 &=3072z^5+3008z^4+z(64z^2-72z+208)+56>0,
\end{align*}
since $64z^2-72z+208>0$ for every real $z$. Equality in \eqref{inv:eq:lower0} therefore occurs exactly at $R=2$.

\smallskip\noindent\emph{Lower endpoint $v=R^2-4$, $R>2$.} In the last case, we obtain
\begin{equation}\label{inv:eq:lowerlarge}
 F(R,R^2-4)-\frac{11}{18}=
 \frac{(R-2)(7R^2+6R-7)}{18(R+1)^2(R+2)}>0.
\end{equation}

The endpoint lemma now shows that $F\ge11/18$ for every admissible pair $(R,v)$. At $R=2,v=0$, formula \eqref{inv:eq:radii} gives $(a,b,c)=(1/3,1/2,1/2)$. At $R=1,v=1$, the admissible sign $u=1$ gives $(a,b,c)=(1/2,1/2,1/3)$; the sign $u=-1$ interchanges the last two circles. The radii are the same in both cases, namely $\{1/3,1/2,1/2\}$, but the order of placement matters. In one configuration the original disks have radii $1/3$ and $1/2$; in the other they both have radius $1/2$. The minimum total area in the theorem is therefore $11\pi/18$.

\section{Third solution: contact-point parametrization}
Our third proof uses contact points as coordinates. In these coordinates, tangency becomes a simple difference condition. This is related to classical chains of tangent circles (see \cite{coxeter68,johnson}). We use Lemma~\ref{lem:global} to check which configurations are allowed.
\subsection{Contact-point coordinates}
Let $O=(0,0)$ be the center of the unit disk, and let $A=(t,0)$ be the center of the first original disk, where $0<t<1$. Since this disk is internally tangent to the unit circle, its radius is $k=1-t$. Here $t$ is a center coordinate, whereas $k$ is a radius; this notation is local to the third solution.  A point on the boundary of the disk centered at $A$ is
\[
P(x)=(t+k\cos x,k\sin x).
\]
The center of a disk of radius $r$ externally tangent to the disk centered at $A$ at $P(x)$ is $P(x)+r(\cos x,\sin x)$.  Its internal tangency to the unit circle means that this center has norm $1-r$.  Solving gives
\[
r=\frac{tk(1-\cos x)}{2-t+t\cos x}.
\]
With $z=\cot(x/2)$ this becomes
\begin{equation}\label{cp:eq:r}
 r(z)=\frac{tk}{k+z^2},\qquad z\in\mathbb R.
\end{equation}
The parameter $z=0$ represents the disk antipodal to the disk centered at $A$, of radius $t$.  The limiting parameter $|z|\to\infty$ represents disks of vanishing radius near the common outer tangency point of the disk centered at $A$.

\begin{lemma}[Tangency as translation]\label{cp:lem:translation}
Two distinct disks from the family \eqref{cp:eq:r}, with parameters $z,w$, are externally tangent if and only if $|z-w|=t$.
\end{lemma}
\begin{proof}
The center of the disk at $z$ is
\[
 C(z)=\left(t+(k+r(z))\frac{z^2-1}{z^2+1},\ (k+r(z))\frac{2z}{z^2+1}\right).
\]
A direct simplification yields
\[
 |C(z)-C(w)|^2-(r(z)+r(w))^2
 =\frac{4k^2((z-w)^2-t^2)}{(k+z^2)(k+w^2)}.
\]
The denominator is positive, so equality is equivalent to $|z-w|=t$.
\end{proof}

Denote by $D_w$ the disk in this family with parameter $w$ and center $C(w)$. By reflection, we may take the second original disk to be $D_z$ with $z\ge0$, and write $B=C(z)$ for its center. The two disks tangent to the disks centered at $A$ and $B$ and to the unit circle are $D_{z-t}$ and $D_{z+t}$. The larger is $D_{z-t}$, because $r(z-t)\ge r(z+t)$.

\subsection{The global maximality condition}
Choosing the larger of the two disks that touch all three circles does not always give the largest disk that fits.  The global admissibility requirement was established in Lemma~\ref{lem:global} using disk radii, and in Lemma~\ref{inv:lem:domain} using inversion parameters.  For this third approach, we express it in contact-point coordinates.  Recall that $A=(t,0)$ is the center and that $k=1-t$ is the radius of the disk centered at $A$.

Apply Lemma~\ref{lem:global} directly to the two original radii $k$ and
$r(z)=tk/(k+z^2)$. In the notation of that lemma, the smaller original
radius must be at least $L$ of the larger one. The disk $D_{z-t}$ is the
larger of the two disks tangent to both originals and the unit circle, so
it is globally largest precisely when this admissibility test holds.

If $k\ge r(z)$, the condition is $r(z)\ge L(k)$. Direct substitution gives
\[
 r(z)-L(k)=
 \frac{t(1-t)(t-z)(t+z)}{(1-t+z^2)(t^2-t+1)}.
\]
Since $z\ge0$, this is equivalent to $z\le t$.

If $r(z)\ge k$, equivalently $z^2\le 2t-1$, then $t\ge1/2$ and
$z\le\sqrt{2t-1}\le t$. The condition is $k\ge L(r(z))$.
Its numerator after clearing the positive denominators factors as
\[
 (k^2+z^2-tz)(k^2+z^2+tz),
\]
up to a strictly positive factor. The second factor is positive, so the
condition is equivalent to $z(t-z)\le k^2$.
Conversely, when $k\ge r(z)$ the inequality $z(t-z)\le k^2$ holds automatically.
Put $g(z)=k^2+z^2-tz$; we need $g(z)\ge0$ whenever $0\le z\le t$ and $z^2\ge 2t-1$.
If $t\le 2/3$, then $k^2\ge t^2/4$, so $g(z)\ge k^2-t^2/4\ge0$ for every $z$.
If $t>2/3$, then $2t-1>t^2/4$, so $z_0=\sqrt{2t-1}>t/2$ and $g$ is increasing on $[t/2,\infty)$;
moreover $g(z_0)=k^2+(2t-1)-tz_0=t(t-z_0)\ge0$ because $z_0\le t$.
Hence $g(z)\ge g(z_0)\ge0$ for all $z\ge z_0$, which is the range in question.
Thus both cases combine into the single admissibility
condition
\begin{equation}\label{cp:eq:admiss}
 0\le z\le t,\qquad z(t-z)\le k^2.
\end{equation}
This is Lemma~\ref{lem:global} expressed in the contact-point coordinates;
no additional geometric maximality argument is needed.

\subsection{Endpoint optimization and equality}
Under \eqref{cp:eq:admiss} the three radii are $k$, $r(z)$, and $r(t-z)$.  Their total area divided by $\pi$ is
\begin{equation}\label{cp:eq:Fz}
 F(t,z)=k^2+\frac{t^2k^2}{(k+z^2)^2}
                +\frac{t^2k^2}{(k+(t-z)^2)^2}.
\end{equation}
Introduce the symmetric parameter
\[
 q=z(t-z),\qquad 0\le q\le q_{\max}:=\min\left\{\frac{t^2}{4},k^2\right\}.
\]
Writing
\[
 H=k^2+k(t^2-2q)+q^2,\qquad S=2k+t^2-2q,
\]
we obtain
\begin{equation}\label{cp:eq:Fq}
 F(t,q)=k^2+t^2k^2\frac{S^2-2H}{H^2}.
\end{equation}

\begin{lemma}[Endpoint principle]\label{cp:lem:endpoint}
For each fixed $t\in(0,1)$, the minimum of $F(t,q)$ on $[0,q_{\max}]$ occurs at an endpoint.
\end{lemma}
\begin{proof}
Differentiation of \eqref{cp:eq:Fq} yields
\[
 F_q=-\frac{4t^2k^2}{H^3}K_t(q),
\]
where
\begin{align*}
 K_t(q)={}&q^3-3(t^2-t+1)q^2\\
 &+(t^4-3t^3+6t^2-6t+3)q-k^3.
\end{align*}
Now
\[
 K_t''(q)=6(q-k^2-t)<0\quad (0\le q\le q_{\max}),
 \qquad K_t'(k^2)=t^3(3-2t)>0.
\]
Since $K_t'$ decreases and $q\le k^2$, we have $K_t'(q)>0$ throughout the interval.  Thus $K_t$ increases strictly, and $F_q$ can change sign only from positive to negative.  Every interior critical point is a strict maximum, not a minimum.
\end{proof}

We now check the three possible endpoints. The factorizations below are exact.

\paragraph{(i) $q=0$.} We have
\begin{equation}\label{cp:eq:end0}
 F(t,0)-\frac{11}{18}=
 \frac{(2t-1)^2 Q_0(t)}{18(t^2-t+1)^2},
\end{equation}
where
\[
 Q_0(t)=9\left(t-\tfrac12\right)^4
 +\tfrac{35}{2}\left(t-\tfrac12\right)^2+\tfrac{33}{16}>0.
\]
Equality holds exactly at $t=1/2$.

\paragraph{(ii) $q=t^2/4$, for $0<t\le2/3$.} We have
\begin{equation}\label{cp:eq:endmid}
 F\left(t,\frac{t^2}{4}\right)-\frac{11}{18}
 =\frac{(3t-2)Q_1(t)}{18(t-2)^4},
\end{equation}
where
\[
 Q_1(t)=6t^5-56t^4+397t^3-618t^2+316t-56.
\]
To verify $Q_1(t)<0$ on $(0,2/3]$, make the substitution $t=2s/[3(1+s)]$ for $s>0$.  Then
\begin{align*}
 -81(1+s)^5Q_1\!\left(\frac{2s}{3(1+s)}\right)
 ={}&1024s^5+3008s^4+4536\\
 &+24s(8s^2-27s+234).
\end{align*}
The quadratic is positive (its discriminant is negative).  Hence \eqref{cp:eq:endmid} is nonnegative, with equality exactly at $t=2/3$.

\paragraph{(iii) $q=k^2$, for $2/3\le t<1$.} Here
\begin{equation}\label{cp:eq:endboundary}
 F(t,k^2)-\frac{11}{18}
 =\frac{(3t-2)(6t^2-20t+7)}{18(t-2)}\ge0.
\end{equation}
Indeed $3t-2\ge0$, while $6t^2-20t+7<0$ and $t-2<0$.  Equality occurs only at $t=2/3$.

By Lemma~\ref{cp:lem:endpoint} and \eqref{cp:eq:end0}--\eqref{cp:eq:endboundary}, $F\ge11/18$.

\paragraph{Both ordered equality configurations.}
The sequential formulation makes it important to distinguish the two
ways equality occurs:
\begin{enumerate}
\item \emph{Two original disks of radius $1/2$.}  Take $t=1/2$ and
$z=0$, so that $k=1-t=1/2$ and
\[
 r(z)=\frac{(1/2)(1/2)}{1/2+0^2}=\frac12,
 \qquad
 r(z-t)=\frac{(1/2)(1/2)}{1/2+(1/2)^2}=\frac13.
\]
The two original radii are $(1/2,1/2)$, and the largest inserted disk
has radius $1/3$.  The maximality criterion holds because
$0\le z=0\le t$ and $z(t-z)=0\le(1-t)^2$.
\item \emph{Original disks of radii $1/3$ and $1/2$.}  Take $t=2/3$
and $z=t/2=1/3$.  Then $k=1/3$ and
\[
 r(z)=r(z-t)=\frac{(2/3)(1/3)}{1/3+(1/3)^2}
 =\frac12.
\]
Thus the original radii are $(1/3,1/2)$, and the largest inserted
disk also has radius $1/2$.  Here $z(t-z)=(1-t)^2=1/9$,
so the maximality criterion holds at equality.
\end{enumerate}
Both cases give the same unordered triple of radii, but they are
distinct \emph{ordered} equality cases because the identity of the
inserted disk differs.  In each case the sum of squared radii is
$1/4+1/4+1/9=11/18$.  Therefore
\begin{equation}\label{cp:eq:answerA}
 {\min\operatorname{Area}=\frac{11\pi}{18},\qquad
 \{r_1,r_2,r_3\}=\left\{\frac13,\frac12,\frac12\right\}.}
\end{equation}

\section{Concluding remarks}
The three proofs look at the same problem in different ways. Descartes' theorem gives an algebraic relation among the radii; inversion turns the picture into equal circles between parallel lines; and contact-point coordinates make tangency easy to describe. Each proof must also check that the inserted disk is the largest disk that fits. Simply choosing the larger of the two disks tangent to the three boundary circles is not enough. If we do not require the inserted disk to touch both original disks, we obtain a different problem.

\end{document}